\input amssymb.sty
\documentclass{amsproc}
\usepackage{graphicx}
\newtheorem{theorem}{Theorem}[section]

\theoremstyle{definition}

\newtheorem{proposition}[theorem]{Proposition}

\theoremstyle{remark}

\numberwithin{equation}{section}

\begin{document}
\title{On the formula of Goulden and Rattan for Kerov polynomials}

\author{Philippe Biane}
\address{CNRS, D\'epartement de Math\'ematiques et Applications,
 \'Ecole Normale Sup\'erieure, 45, rue d'Ulm 75005 paris, FRANCE
}
\email{Philippe.Biane@ens.fr}
\subjclass{Primary ; Secondary }
\date{}

\begin{abstract}
We give a simple proof of an explicit formula for Kerov polynomials.
This explicit formula is closely related to a recent formula of Goulden and
Rattan.
\end{abstract}

\maketitle

\section{Kerov polynomials}
Kerov polynomials are universal polynomials which express the characters of
symmetric groups evaluated on cycles, in terms of quantities known as the free
cumulants of a Young diagram. We now explain these notions.\par
Let $\lambda=\lambda_1\geq \lambda_2\geq\ldots$ be a Young diagram, to
 which we associate a piecewise affine function 
 $\omega:\mathbb R\to \mathbb R$, with slopes $\pm1$,
 such that $\omega(x)=|x|$ for $|x|$ large enough, as in Fig. 1 below, which
 corresponds to the partition $8=4+3+1$.
  We can encode the Young diagram using the local minima
and local maxima of the function $\omega$, denoted by $x_1,\ldots, x_m$ and
$y_1,\ldots, y_{m-1}$ respectively, which form two interlacing sequences of
integers. These are (-3,-1,2,4) and (-2,1,3) respectively in the picture. 
$$\includegraphics{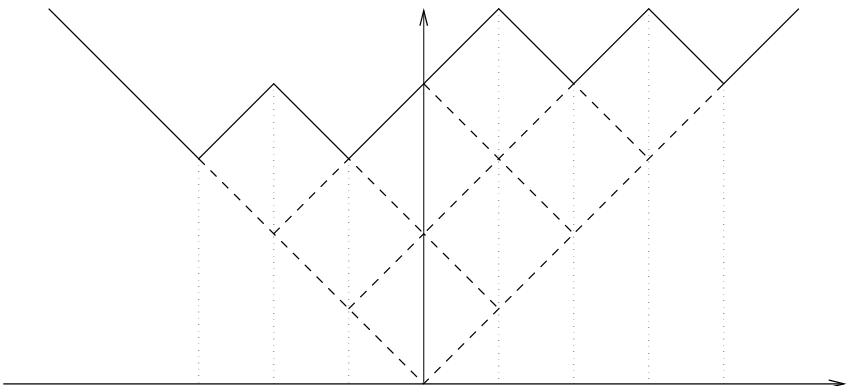}$$
  $$\qquad\qquad\ x_1\quad\,  y_1\quad 
x_2\qquad\qquad  y_2\quad\, x_3\quad y_3\quad x_4\qquad$$
$$Fig.1$$
To the Young diagram we associate the rational fraction
$$G_{\lambda}(z)={\prod_{i=1}^{m-1}(z-y_i)\over
\prod_{i=1}^{m}(z-x_i)}
$$
and the formal power series $K_{\lambda}$,  inverse to $G$ for composition, 
$$K_{\lambda}(z)=G_{\lambda}^{(-1)}(z)=z^{-1}+\sum_{k=2}^\infty
R_k(\lambda)z^{k-1}.
$$
 The quantities $R_k(\lambda);k=2,3,\ldots$ are called the free cumulants of 
  the diagram $\lambda$. Note that $R_1(\lambda)=0$ for any Young diagram, so we
  do not include it
   in the series of free cumulants. These quantities arise in the
  asymptotic study of representations of symmetric groups, see 
  \cite{B1}.
  It turns out that there exists universal polynomials
   $\Sigma_2,\Sigma_3\ldots$ in the variables $R_2,R_3,\ldots$ such that for any
   Young diagram $\lambda$ the normalized character $\chi_{\lambda}$ evaluated
   on a 
   cycle of length $k$ is given by 
   $$(n)_k\chi_{\lambda}(c_k)=\Sigma_k(R_2(\lambda),R_3(\lambda),\ldots).$$
   The remarkable fact here is that these polynomials do not depend on the size
   of the symmetric group.
    We list the  first such polynomials below.
 $$
\begin{array}{l}
\Sigma_1=R_2\\
\Sigma_2=R_3\\  
\Sigma_3=R_4+R_2\\ 
 \Sigma_4=R_5 +5R_3\\ 
  \Sigma_5=R_6+15R_4+5R_2^2+8R_2\\   
\end{array}
$$
We refer to \cite{B2} and \cite{GR}
for more information about results and conjectures
on the coefficients of these polynomials. We take from \cite{B2}, section 5, the
following expression for Kerov polynomials. Here $[z^{-k}]\,f(z)$
denotes the coefficient of $z^{-k}$ (the residue if $k=1$) of a Laurent series
$f(z)$.
\begin{proposition} Consider the formal power series
  $$H(z)=z-\sum_{j=2}^\infty B_jz^{1-j}.$$
Define $$\Sigma_{k}=-\frac{1}{k}\,[z^{-1}]\,
H(z)H(z-1)\ldots H(z-k+1)\qquad (1)$$
and
$$R_{k+1}=
-{1\over k}\,[z^{-1}]\,
H(z)^k$$
then the expression of $\Sigma_k$ in terms of the $R_k\,'$s is given by Kerov's
polynomials.
\end{proposition}
Recently I. P. Goulden and A. Rattan \cite{GR}
 have given an explicit expression for Kerov
polynomials, from which they have deduced a certain number of positivity
properties of the coefficients of these polynomials. Their proof uses the
Lagrange inversion formula. In the next section we use the invariance of residue
under change of variables to 
derive in a simple way a closely related 
formula, and show how to recover Goulden
and Rattan's formula.
\section{Explicit expression for Kerov polynomials}
We use the notations of Proposition 1.1 above.
Let us introduce the power series
$$L(z)=z+\sum_{j=2}^\infty R_jz^{1-j}.$$
One has $H\circ L(z)=z$, by Lagrange inversion formula.
We use the invariance of the residue under change of variables, namely if
$u,f$ are Laurent series, and $u$ is invertible for composition, then
$$[z^{-1}]\,f(z)=[\zeta^{-1}]\,u'(\zeta)f\circ u(\zeta).$$
 Using Taylor's formula as well as the change of variables $z=L(\zeta)$
in the residue, one gets from (1)

$$
\begin{array}{rcl}
\Sigma_{k}&=&-{1\over k}\,[z^{-1}]\,\prod_{j=0}^{k-1}
\left(\sum_{r=0}^{\infty}\frac{(-j)^r}{r!}H^{(r)}(z)\right)
\\
&=&-{1\over k}\,[\zeta^{-1}]\,L'(\zeta)\prod_{j=0}^{k-1}
\left(\sum_{r=0}^{\infty}\frac{(-j)^r}{r!}H^{(r)}\circ L(\zeta)\right).
\end{array}
$$
Using $H'\circ L(\zeta)=\frac{1}{L'(\zeta)}$ one gets
$H^{(r)}\circ
L(\zeta)=\left(\frac{1}{L'(\zeta)}\frac{d}{d\zeta}\right)^{r-1}
\frac{1}{L'(\zeta)}
$
therefore
$$\Sigma_{k}
=-{1\over k}\,[\zeta^{-1}]\,L'(\zeta)\prod_{j=0}^{k-1}
\left(\zeta+\sum_{r=1}^{\infty}\frac{(-j)^r}{r!}
\left(\frac{1}{L'(\zeta)}\frac{d}{d\zeta}\right)^{r-1}
\frac{1}{L'(\zeta)}\right).
$$
Putting $F(\zeta)=\frac{1}{L'(\zeta)}$ 
we obtain the following Proposition.
\begin{proposition}
Let $$F(\zeta)=\frac{1}{L'(\zeta)}=\frac{1}{1-\sum_{k=2}^\infty
(k-1)R_k\zeta^{-k}}$$ then Kerov's polynomials are given by the following
expression
$$
\Sigma_{k}=-{1\over k}\,[\zeta^{-1}]\,\frac{1}{F(\zeta)}\prod_{j=0}^{k-1}
\left(\zeta+\sum_{r=1}^{\infty}\frac{(-j)^r}{r!}
\left(F(\zeta)\frac{d}{d\zeta}\right)^{r-1}
F(\zeta)\right).\qquad (2)
$$

\end{proposition}
\section {The formula of Goulden and Rattan}
Goulden and Rattan give various equivalent formulas for $\Sigma_k$.
They introduce the series $C(\zeta)=F(\zeta^{-1})$, and
define polynomials
$$P_m(z)=-\frac{1}{m!}C(z)(D+(m-2)I)\left[C(z)\ldots 
(D+I)\left[C(z)DC(z)\right]\ldots\right]$$ 
 where $D=z\frac{d}{dz}$.
 The generating series form of their formula now reads

$$\Sigma_k=-\frac{1}{k}[z^{k+1}]\frac{1}{C(z)}\prod_{j=1}^{k-1}(1+\sum_{i=1}^\infty
j^iP_i(z)z^i).\qquad (3)$$
We  recover this formula using $(2)$. For this we
factor out $\zeta^k$ in the expression in the rhs of $(2)$, and use the
change of variable $z=\zeta^{-1}$. This gives 
$$
\begin{array}{rcl}
\Sigma_{k}&=&-{1\over k}\,[\zeta^{-1}]\,\frac{1}{F(\zeta)}\prod_{j=0}^{k-1}
\left(\zeta+\sum_{r=1}^{\infty}\frac{(-j)^r}{r!}
\left(F(\zeta)\frac{d}{d\zeta}\right)^{r-1}
F(\zeta)\right)\\
&=&-{1\over k}\,[\zeta^{-1-k}]
\,\frac{1}{F(\zeta)}\prod_{j=1}^{k-1}
\left(1+\sum_{r=1}^{\infty}\frac{(-j)^r}{r!}\zeta^{-1}
\left(F(\zeta)\frac{d}{d\zeta}\right)^{r-1}
F(\zeta)\right)\\
&=&-{1\over k}\,[z^{k+1}]\,\frac{1}{C(z)}\prod_{j=1}^{k-1}
\left(1+\sum_{r=1}^{\infty}\frac{(-j)^r}{r!}z
\left(-C(z)z^2\frac{d}{dz}\right)^{r-1}
C(z)\right)\\
&=&-{1\over k}\,[z^{k+1}]\,\frac{1}{C(z)}\prod_{j=1}^{k-1}
\left(1-\sum_{r=1}^{\infty}\frac{j^rz}{r!}(C(z)zD)^{r-1}C(z)
\right).\\
\end {array}
$$
Now remark that $z^{-j}\circ D\circ z^j=D+jI$ to get
$$\begin{array}{rcl}
z(C(z)zD)^{r-1}C(z)&=&z^rC(z)(D+(r-2)I)\left[C(z)\ldots 
(D+I)\left[C(z)DC(z)\right]\ldots\right]\\
&=&-r!P_r(z).
\end{array}$$
\bibliographystyle{amsalpha}

\end{document}